\documentclass[xypic,amssymb,rawfonts,amsthm,amsmath,amsfonts,syntonly,11pt]{amsart}

\usepackage{amssymb}

\newcounter{thecounter}
\numberwithin{thecounter}{section}
\newtheorem{lemma}[thecounter]{Lemma}
\newtheorem{prelemma}[thecounter]{Pre-Lemma}

\newtheorem{thm}[thecounter]{Theorem}

\theoremstyle{definition}
\newtheorem{rec}[thecounter]{Recollection}

\newtheorem{rem}[thecounter]{Remark}

\numberwithin{equation}{section}

\input xypic
\xyoption{all}
\CompileMatrices

\newcommand{\sd}{\operatorname{sd}}

\newcommand{\sta}{\operatorname{\bf star}}
\newcommand{\link}{ \operatorname{\bf link}}

\newcommand{\Rmod}{R\mbox{-}\operatorname{mod}}

\newcommand{\fF}{{\mathfrak F}}

\newcommand{\cC}{{\mathcal{C}}}
\newcommand{\calD}{{\mathcal{D}}}
\newcommand{\cX}{{\mathcal{X}}}

\newcommand{\bO}{{\mathbf O}}

\newcommand{\op}{\operatorname{op}}

\newcommand{\GL}{\operatorname{GL}}

\newcommand{\SL}{\operatorname{SL}}

\newcommand{\F}{{\mathbf {F}}}
\newcommand{\Z}{{\mathbf {Z}}}

\newcommand{\cE}{{\mathcal{E}}}
\newcommand{\A}{{\mathbf A}}

\newcommand{\beq}{\begin{eqnarray*}}
\newcommand{\eeq}{\end{eqnarray*}}
\newcommand{\pil}[1]{\stackrel{#1}{\rightarrow}}

\newcommand{\tuborg}{\left\{\begin{array}{ll}}
\newcommand{\sluttuborg}{\end{array}\right.}

\newfont{\bm}{msbm10}
\newcommand{\semi}{{\mbox{{\bm \symbol{111}}}}}

\newdir{ >}{{}*!/-5pt/\dir{>}} %\ar@{ >->} gives a monomorphism arrow

\newcommand{\cB}{{\mathcal B} }

\newcommand{\cZ}{{\mathcal Z} }
\newcommand{\cA}{{\mathcal A} }
\newcommand{\cS}{{\mathcal S} }
\newcommand{\cI}{{\mathcal I} }
\newcommand{\cCe}{{\mathcal Ce} }
\newcommand{\cBCe}{{{\mathcal B}{\mathcal C}{e}} }
\newcommand{\Id}{\operatorname{Id}}

\begin{document}
\title[Propagating sharp decompositions]
      {Propagating sharp group homology decompositions}
%\date{\today}

\author[J.~Grodal]{Jesper Grodal}
\author[S.~D.~Smith]{Stephen D.~Smith}
\thanks{The first author was partially supported
        by NSF grant DMS-0104318.
	The second author was partially supported
	by NSA grant MDA 904-01-1-0045.  }
\subjclass{Primary: 20J06; Secondary: 20J05, 55R35, 55P91}
%\thanks{\bf To Advances in Math: Please send editorial correspondence to Jesper Grodal, email jg@math.uchicago.edu (fax (+1) (773) 702 9787 ; regular mail: Dept. of Mathematics, University of Chicago, Chicago, IL 60637, USA)}
\keywords{ cohomology decomposition, induction theory, subgroup complex}
\address{Department of Mathematics, University of Chicago, Chicago, IL 60637, USA}
\email{jg@math.uchicago.edu}

\address{Department of Mathematics (m/c 249),
         University of Illinois at Chicago,
	 851 S. Morgan, 
         Chicago, IL 60607-7045, USA}
\email{smiths@math.uic.edu}

\begin{abstract}
A collection $\cC$ of subgroups of a finite group $G$
can give rise to three different standard
formulas for the cohomology of $G$ in terms of either
the subgroups in $\cC$ or their centralizers or their normalizers.
We give a short but systematic study
of the relationship among such formulas
for nine standard collections $\cC$ of $p$-subgroups,
obtaining some new formulas in the process.
To do this, we exhibit some sufficient conditions on the poset $\cC$
which imply comparison results.
\end{abstract}

\maketitle

\section{Introduction and statement of the two theorems}
            \label{sec:introthms}
Induction theory,
for representations of a group $G$ over a commutative ring $R$,
is about calculating the value $\fF(G)$ of a functor $\fF$
from the opposite orbit category $\bO(G)^{\op}$ of $G$ to $R$-modules,
in terms of the values $\fF(H)$ for various subgroups $H$ of $G$.
From ``topological induction theory'',
which investigates homology decompositions, one gets,
for a fixed collection $\cC$ of subgroups of $G$,
three different potential induction formulas:
by considering the subgroups $H$ themselves, their normalizers,
or their centralizers.
Furthermore, restricting attention now to collections of (non-trivial)
$p$-subgroups of a finite group $G$, for a fixed prime $p$,
there are at least nine different collections commonly studied,
which thus give us some $27$ potential induction formulas.

The goal of this paper is to give a brief but systematic study
of the interrelationship among these $27$ different potential formulas.
Namely, we investigate when the validity of one formula
implies the validity of other formulas, for geometric reasons,
reasons which do not depend on the specific functor $\fF$,
or only depend on $\fF$ in a mild way---this
is what we mean by ``propagation'' in the title of the paper.
We proceed by studying homotopy properties
of certain associated poset spaces,
extending an approach of Dwyer \cite{dwyer97,dwyer98sharp};
see also \cite{grodal02}.
Our proofs basically consist of three short lemmas
in Section \ref{sec:lemmas},
which each state assumptions on the poset $\cC$
which guarantee comparison results among the formulas.
As a consequence we are able to settle,
either in the positive or the negative,
which of these $27$ formulas give rise to 
``sharp homology decompositions'' \cite{dwyer98sharp}---recovering 
many previous results, as well as closing the gaps
in the existing literature.

\medskip

For some collection of subgroups $\cC$
(always assumed closed under $G$-conjugation),
let $\bO_\cC$ denote the full subcategory of the orbit category $\bO(G)$
with objects $G/H$,
where we take only those $H$ which are members of $\cC$.
The most na\"ive approach to induction theory is to ask if the functor 
  $$ F^\beta: \bO_\cC^{\op} \to \Rmod,
        \mbox{ given by } H \mapsto \fF(H)$$
satisfies that
  $$ \fF(G) \pil{\cong} \textstyle{\lim_{\bO_\cC}^0}\ F^\beta , $$
i.e., that $\fF$ is {\em (subgroup) $\cC$-computable}
in the language of \cite[p.293]{dress75}.
A more ambitious demand is to require that also
  $$ \textstyle{\lim_{\bO_\cC}^i}\ F^\beta = 0 \mbox{ for } i > 0 , $$
where $\lim^i$ denotes the $i$th higher derived functor
of the inverse limit functor;
in which case $\fF$ is called {\em subgroup $\cC$-acyclic}
\cite[Def.~8.5]{grodal02}.
See also \cite{grodal02} for motivation.

One can ask the same for two corresponding functors
obtained via normalizers and centralizers:
  $$ F^\delta: ((\sd \cC)/G)^{\op} \to \Rmod, \mbox{ given by }
     ( H_0 < \cdots < H_k ) \mapsto
                 \fF \bigl( N_G(H_0) \cap \cdots \cap N_G(H_k) \bigr) $$
  $$ \mbox{ and }
     F^\alpha: \A_\cC \to \Rmod \mbox{ given by }
                                 H \mapsto \fF \bigl( C_G(H) \bigr) . $$
Here $(\sd \cC)/G$ is the orbit simplex category, with objects
the $G$-conjugacy classes of strict chains of subgroups in $\cC$,
and morphisms the refinement of chains;
and $\A_\cC$ is the conjugacy category,
with objects the subgroups in $\cC$
and morphisms the homomorphisms between subgroups
induced by conjugation in $G$.

Independently of whether a collection $\cC$ gives induction formulas
in the above senses of computability or acyclicity,
it is of interest to study the relationship among the limits 
$\lim^*_{\bO_\cC} F^\beta$, $\lim^*_{ (\sd \cC)/G} F^\delta$,
and $\lim^*_{\A_\cC} F^\alpha$,
as well as to examine how these limits change
when the collection $\cC$ is varied.
Aspects of this question have already been much studied
because of applications in homotopy theory,
where these higher limits often play an important role;
see e.g., \cite{dwyer-henn01} and \cite{grodal02}, and their references.

Using the definitions it is possible
(see e.g., \cite[\S 1.2]{dwyer97}, \cite[2.8]{grodal02}
or \cite[\S 10]{dwyer-henn01})
to reformulate these higher limits in terms of Bredon cohomology,
which makes dealing with them more amenable to homotopy theory.
More precisely, let $E\bO_\cC$ be the poset category
with objects $(G/H,gH)$, where $H \in \cC$,
and morphisms from $(G/H,x)$ to $(G/H',x')$ given by the
$G$-maps $G/H \to G/H'$ sending $x$ to $x'$;
and let $E\A_\cC$ be the category
with objects the monomorphisms $i: H \to G$,
and morphisms from $i$ to $i'$ given by
the homomorphisms $\varphi: H \to H'$ such that $i = \varphi i'$.
These spaces admit natural $G$-actions.
With this notation we have the following models for the higher limits:
  $$ \lim_{\bO_\cC}{}^* F^\beta = H^*_G(|E\bO_\cC|;\fF), \mbox{  } 
     \lim_{ (\sd \cC)/G} \hspace{-8pt}{}^* F^\delta = H^*_G(|\cC|;\fF),
          \mbox{ and }
     \lim_{\A_\cC}{}^* F^\alpha = H^*_G(|E\A_\cC|;\fF) . $$
Here $H^*_G(X;\fF)$ the denotes Bredon cohomology of a space $X$
with values in the generic coefficient system $\fF$---in other words,
$H^*_G(X;\fF)$ is the homology of the canonical cochain complex
which in degree $n$ is given
by $\prod_{[\sigma] \in X_n/G}\ \fF(G_\sigma)$,
where $G_\sigma$ is the stabilizer of the $n$-simplex $\sigma$.
(Different, smaller, models are found in \cite{grodal02}
under assumptions on $\cC$.)

Note that we have natural comparison functors
  $$ E\bO_\cC \to \cC \leftarrow E\A_\cC $$
given by $(G/Q,x) \mapsto G_x$ and $(i: Q \to G) \mapsto i(Q)$.
These functors are in fact equivalences of categories,
with inverses being given
by $Q \mapsto (G/Q,eQ)$ and $Q \mapsto (incl: Q \to G)$,
and so the three categories have homotopy equivalent nerves.
However these inverses generally do not respect the $G$-action;
and the nerves are in general not $G$-homotopy equivalent.

More generally, for a subcollection $\cC'$ of $\cC$
we can consider the diagram of $G$-spaces and $G$-maps
  $$ \xymatrix{ |E\bO_{\cC'}| \ar[r] \ar[d] & |{\cC'}| \ar[d]
                                     & \ar[l] {|E\A_{\cC'}|} \ar[d]\\
         |E\bO_{\cC}| \ar[r] & |{\cC}| & \ar[l] {|E\A_{\cC}|}} $$

\noindent
If a map in this diagram is a $G$-homotopy equivalence,
then it induces an isomorphism between the corresponding higher limits,
by applying $H^*_G(-;\fF)$.
In some applications, $\fF$ has the further structure
of a cohomological Mackey functor \cite{yoshida83}---that is,
a Mackey functor with the standard property \cite[II.5.3]{AM94}
of group cohomology functors $H^n(-;M)$
that restriction to a subgroup $H$ followed by transfer back to $G$
is just multiplication by the index $|G:H|$.
Then if $|G:H|$ is also assumed invertible in $R$, we only need
an $H$-equivalence, to get an induced isomorphism above.
In particular when $R = \Z_{(p)}$ it is enough
to have an $S$-equivalence, for $S$ a Sylow $p$-subgroup of $G$.

Our main theorem \ref{equivalencethm} below
concerns 9 particular collections $\cC$ of $p$-subgroups,
whose definitions we review after the statement of the theorem.
(These collections are interesting
because of their known acyclicity and ``sharpness'' properties---see
Theorem \ref{sharpnessthm}.)
The theorem says exactly when the maps
in the diagram of the previous paragraph for these $\cC$
are either $G$-equivalences,
or $S$-equivalences for $S$ a Sylow $p$-subgroup of $G$.

\begin{thm} \label{equivalencethm} Fix an arbitrary prime $p$,
and let $G$ denote a finite group, with Sylow $p$-subgroup $S$.
Then we have the following table:
$$\xymatrix@-15pt{ 
%  *[F-:<4pt>]{ \mbox{ space for }  \cC \rightarrow  \mbox{}} 
        & \calD & \cBCe & \cCe & \cB & \cI & \cS & \cA & \cZ &  \cE   \\
 {|E\bO_\cC| : } & \bullet \ar@{..}[d]
                 & \bullet  \ar@{--}[d]  \ar@{-}[r]
        & \bullet  \ar@{--}[d] & \bullet  \ar@{--}[d]  \ar@{-}[r]
        &  \bullet  \ar@{--}[d]  \ar@{-}[r]
        & \bullet \ar@{--}[d]  \ar@{..}[r]
        & \bullet  \ar@{-}[r] \ar@{..}[d] & \bullet \ar@{..}[d]
        & \bullet  \ar@{..}[d]\\
{|\cC| : }  & \bullet \ar@{..}[d] & \bullet  \ar@{-}[r] \ar@{..}[d]
        & \bullet \ar@{..}[d] & \bullet \ar@{-}[r] \ar@{..}[d]
        &  \bullet    \ar@{-}[r] \ar@{..}[d]
        & \bullet \ar@{--}[d]  \ar@{-}[r]
        & \bullet  \ar@{--}[d]  \ar@{-}[r] & \bullet  \ar@{--}[d]
        & \bullet  \ar@{--}[d] \\
{|E\A_\cC| : }  & \bullet & \bullet \ar@{..}[r]
       & \bullet & \bullet \ar@{..}[r]
       & \bullet \ar@{..}[r]
       & \bullet  \ar@{-}[r] & \bullet  \ar@{-}[r]
       & \bullet & \bullet} $$
Here a node denotes the space indicated by its row,
for the collection $\cC$ indicated by its column. 
A dotted line denotes a homotopy equivalence;
a dashed line denotes an $S$-homotopy equivalence;
and a solid line denotes a $G$-homotopy equivalence.

These results are best-possible in the sense
that if a certain line is not present,
then (for any prime $p$) there exists a finite group $G$
for which that kind of equivalence does not hold.
\end{thm}

\noindent
We will later in Remark~\ref{histrem1} comment on the history
behind various parts of this result---in particular,
the horizontal lines in the middle row
correspond to classical equivalence theorems in the literature.

\medskip

{\em Notation for collections of nontrivial $p$-subgroups:}
The collection $\cS = \cS_p(G)$
is the collection of all non-trivial $p$-subgroups \cite{brown75}.
Next $\cB = \cB_p(G)$ is the subcollection of $\cS$
given by all non-identity {\em $p$-radical subgroups} \cite{bouc84}
i.e., all non-trivial $p$-subgroups $Q$
such that $O_p \bigr( N_G(Q)/Q \bigl)\ = 1$.
Also $\cCe = \cCe_p(G)$
is the subcollection of $\cS$ of {\em $p$-centric subgroups}
\cite{dwyer97}, i.e., $p$-subgroups $Q$
such that $Z(Q)$ is a Sylow $p$-subgroup in $C_G(Q)$.
Then we let $\cBCe = \cB \cap \cCe$ denote
the collection of nontrivial $p$-central and $p$-radical subgroups.
Further $\calD = \calD_p(G)$ is the subcollection of $\cBCe$
given by {\em principal} $p$-radical subgroups \cite{grodal02},
i.e., the subgroups $Q$ in $\cCe$
such that $O_p \bigl( N_G(Q) / Q C_G(Q) \bigr) = 1$.
We let $\cI$ denote the subcollection of $\cS$
given by the {\em Sylow intersections}: all non-trivial subgroups
which are intersections of a set of Sylow $p$-subgroups in $G$.
Next $\cA = \cA_p(G)$ is the subcollection of $\cS$ of all nontrivial
{\em elementary abelian $p$-subgroups} \cite{quillen78}.
Then $\cZ = \cZ_p(G)$ is the subcollection of $\cA$
in \cite[Sec.~6.6]{benson98}, given by subgroups $V$
such that $\Omega_1 O_p Z \bigl( C_G(V) \bigr) = V$
(the former expression denotes the elements of order dividing $p$
in the center of the centralizer of $V$).
Finally $\cE = \cE_p(G)$ \cite[Sec.~3]{benson94}
is the smallest subcollection of $\cA$
which contains the conjugates of the subgroups of order $p$
in the center of a Sylow $p$-subgroup of $G$,
and is closed under taking products of commuting members.

\medskip

We recall that a collection $\cC$ is said
to give rise to a {\em sharp subgroup (co)homology decomposition},
or to be {\em subgroup sharp} for short,
if  $H^n(-;\F_p)$ is subgroup $\cC$-acyclic for all $n \geq 0$, i.e., 
if $\lim_{G/Q \in \bO_\cC}^i H^n(Q;\F_p) = 0$ for $i > 0$
and $H^n(G;\F_p) \pil{\cong} \lim_{G/Q \in \bO_\cC}^0 H^n(Q;\F_p)$
for all $n \geq 0$.
{\em Normalizer} and {\em centralizer sharpness} are defined similarly;
we refer to e.g.\ \cite{dwyer98sharp}, \cite[\S 8]{dwyer-henn01},
and \cite[\S 9]{grodal02} for background and motivation.

Combining Theorem~\ref{equivalencethm} with known results,
we are able to complete the following table,
showing which types of sharpness hold
for each of the collections studied above---obtaining
alternative proofs of sharpness for many of them in the process.
Here the row names
are abbreviations for subgroup, normalizer, and centralizer sharpness;
a ``$y$'' in the table
indicates that sharpness holds for all finite groups $G$ and all primes $p$;
while ``$n$'' indicates that for any prime $p$,
sharpness fails for some $G$.
(If a positive result was previously known, we have given a reference
to the place where it first seemed to be stated in the literature.)
\begin{thm} \label{sharpnessthm}
{\small
$$\begin{array}{l||l|l|l|l|l|l|l|l|l}
        & \calD & \cBCe & \cCe & \cB & \cI & \cS & \cA & \cZ &  \cE   \\
	    \hline \hline
\mbox{s:}               &   y 
        &    y                       &   y                       
        &    y                       & y 
        & y                        &  n
        &   n   &   n  \\
                 &  \mbox{\cite{grodal02}}  
        &      \mbox{\cite{dwyer97}} &     \mbox{\cite{dwyer97}}
        &      \mbox{\cite{dwyer97}} &   \mbox{\cite{notbohm01}}
        &   \mbox{\cite{dwyer97}}  &   
        &       &      \\
	\hline
\mbox{n:}            &   y  
        &    y                        &   y                        
        &   y                      &  y                         
        & y                       &  y                     
        &   y                          &   y   \\
                              &      \mbox{\cite{grodal02}}
        &      \mbox{\cite{grodal02}} &      \mbox{\cite{grodal02}}
        &     \mbox{\cite{webb91}} &    \mbox{\cite{notbohm01}}
        &   \mbox{\cite{webb91}}  &    \mbox{\cite{webb91}}
        &       \mbox{\cite{benson98}} &       \\
	\hline
\mbox{c:}            &   n 
        &    n  &   n
        &    n  &   n 
        &  y                       &  y                       
        &   y  & y  \mbox{ \cite{JM92},}                \\
        &     
        &       &    
        &       &     
        &    \mbox{\cite{dwyer97}} &        \mbox{\cite{JM92}}
        &      &   \mbox{\cite{benson94}} 
\end{array}
$$
}
\end{thm}

In particular, observe that the new positive results here
are that the collection $\cZ$ is centralizer sharp,
and that the collection $\cE$ is normalizer sharp.
These are of some interest,
since both $\cZ$ and $\cE$ are subcollections---often proper---of the
more common collection $\cA$
of all nontrivial elementary abelian $p$-subgroups;
for example, the normalizer sharpness of $\cE$ is applied
to a number of sporadic simple groups at $p = 2$ in \cite{BS03}.

 A smaller version of the table in \ref{sharpnessthm} was given
 by Dwyer at his talk at the $1996$ AMS Summer Research Institute
 on Representations and Cohomology;
 and the present note arose
 as an attempt to go back and complete and extend that table.  

\section{Three lemmas and a pre-lemma}   \label{sec:lemmas}

We start with some recollections and notation.
Recall  that a map is a $G$-homotopy equivalence
if and only if the induced map on $H$-fixed points is an (ordinary)
homotopy equivalence for all subgroups $H \leq G$
(see \cite[6.4.2]{benson98}).
We say that a poset $\cX$ is contractible
if its nerve $|\cX|$ is contractible.
For a poset $\cX$ and $x \in \cX$, $\cX_{\leq x}$
is the subposet of elements less than or equal to $x$,
with $\cX_{\geq x}$, $\cX_{< x}$, and $\cX_{> x}$ defined analogously.
Further define
$\sta_\cX(x) = \{ y \in \cX | x \leq y \mbox{ or } y \leq x \}$
and $\link_\cX(x) = \{ y \in \cX | x < y \mbox{ or } y < x \}$;
the nerves of these posets are of course the star and link
of the vertex $x$ in the nerve $|\cX|$ of $\cX$.
Note that
  $$ \link_\cX(x) = \cX_{<x}\ \star\ \cX_{>x} \mbox{ and }
     \sta_\cX(x) = \cX_{< x} \star\ x \star\ \cX_{>x} ,
                                                     \eqno{(\star)} $$
where $\star$ denotes the join of posets,
obtained by taking the disjoint union as sets
and imposing the additional order relation
that all elements in left poset are smaller
than the ones in the right poset.
By \cite[Prop.~1.9]{quillen78}, on nerves the join of posets
produces the join of spaces, which we also denote by $\star$. 
(For example we see using remarks above
that we have a $G_x$-homeomorphism
$ |\sta_\cX(x)| \approx |x|\ \star\ |\link_\cX(x)|$.)

We of course have that
  $$ \cC^H = \{ Q \in \cC\ |\ H \leq N_G(Q) \}  . $$
The elementary but fundamental observation
which is needed for the lemmas of this section
is that the functors $E\bO_\cC \to \cC \leftarrow E\A$
described in the introduction induce equivalences of categories
  $$ \xymatrix@-10pt{ E\bO_{\cC}^H  \ar[r]
       & \cC_{\geq H}   & \text{ and }  \\
       &    \cC_{\leq C_G(H)}
                      & E\A_\cC^H \ar[l] }  \eqno{(\dagger)}$$ 
which are natural in the variable $\cC$.
This can be used to examine the more precise failure
of the functors  $E\bO_\cC \to \cC \leftarrow E\A_\cC$
to induce $G$-homotopy equivalences.

\begin{rec} \label{quillenobs}
The most fundamental fact in the homotopy theory of categories
is the observation that a natural transformation between two functors
induces a homotopy between their nerves
(see e.g., \cite[1.3]{quillen78}).
In particular if $F$ is an endomorphism of a poset $\cX$
that is order-related to the identity functor,
i.e., if either $F \leq \Id_\cX$
(that is, $F(x) \leq x$ for all $x \in \cX$) or $F \geq \Id_\cX$,
then we can construct a natural transformation from $F$ to $\Id_\cX$
or vice versa.
Hence $F$ induces a homotopy deformation retraction
of the nerve of $\cX$ onto the nerve of any poset $\cX'$
such that $F(\cX) \subseteq \cX' \subseteq \cX$.
If $\cX'$ is a $G$-subposet of a $G$-poset $\cX$
and $F$ is $G$-equivariant,
then the retraction will be a $G$-homotopy deformation retraction,
and $|\cX'|$ is $G$-homotopy equivalent to $|\cX|$.
In particular a poset with a unique largest or smallest element $x$
is contractible, since we can take $F$ to be the 
endomorphism sending everything to $x$.
\end{rec}

While the above technique is obviously useful
to get results about $\cC$ from functors $F$ on $\cC$,
it can also be propagated
to uses on $E\bO_\cC$ and $E\A_\cC$ via $(\dagger)$.

\begin{lemma} \label{lemma1}
Let $\cC$ be a collection of subgroups of a discrete group $G$.
Suppose that $F$ is a $G$-equivariant poset endomorphism of $\cC$
satisfying either $F \geq \Id_{\cC}$ or $F \leq \Id_{\cC}$.
Set $\cC' = F(\cC)$.

\begin{enumerate}

\item 
Assume that the case $F \geq \Id_\cC$ of the hypothesis holds. \\
Then the inclusion $E\bO_{\cC'} \to E\bO_\cC$
induces a $G$-homotopy equivalence on nerves.

\item 
Assume for all $P \in \cC$
that $C_G(P) \leq C_G \bigl( F(P) \bigr)$. \\
Then the inclusion $E\A_{\cC'} \to E\A_{\cC}$
induces a $G$-homotopy equivalence on nerves. \\
(Notice that the further hypothesis holds
in the case $F \leq \Id_\cC$ of the hypothesis.)

\item The inclusion $\cC' \to \cC$
induces a $G$-homotopy equivalence on nerves.
\end{enumerate}

\end{lemma}

\begin{proof}
First note that (3) follows directly from \ref{quillenobs}.

To see (1), we want to see that $E\bO_{\cC'}^H \to E\bO_\cC^H$
induces a homotopy equivalence on nerves, for all subgroups $H \leq G$.
For this, observe by ($\dagger$) that this map identifies
with the inclusion $\cC'_{\geq H} \to \cC_{\geq H}$,
up to equivalence of categories.
The assumption $F \geq \Id_\cC$ guarantees that $F$ takes
$\cC_{\geq H}$ into $\cC_{\geq H}$,
and so induces a homotopy deformation retraction
of $|\cC_{\geq H}|$ onto $|\cC'_{\geq H}|$ by \ref{quillenobs}.

Finally, (2) follows by a parallel argument,
since now the assumptions on $F$ guarantee
that $F$ induces a homotopy deformation retraction of
the nerve of $\cC_{\leq C_G(H)}$ onto the nerve of $\cC'_{\leq C_G(H)}$:
since if $Q \leq C_G(H)$, then $H \leq C_G(Q) \leq C_G(F(Q))$,
so that $F(Q) \leq C_G(H)$.
\end{proof}

Next, we work towards a lemma that says when we can remove subgroups,
one conjugacy class at a time,
from a collection, while preserving $G$-homotopy type---again
we propagate standard methods from $\cC$ to $E\bO_\cC$ and $E\A_\cC$.
First a recollection and a pre-lemma.
\begin{rec} \label{pushoutrec}
Let $\cX$ be a $G$-poset, and $\cX'$ the subposet of $\cX$
obtained by removing from $\cX$ the $G$-conjugates of an element $x$.
Then we have the following (homotopy) pushout-square of $G$-spaces:
  $$ \xymatrix{ G \times_{G_x} |\link_\cX(x)| \ar[r] \ar[d]
                                              & |\cX'|  \ar[d] \\
              G \times_{G_x} |\sta_\cX(x)| \ar[r] & |\cX| } . $$
To see that we have a pushout square of spaces,
just observe that every simplex in $|\cX|$ which is not in $|\cX'|$
lies in the $G$-orbit of $|\sta_\cX(x)|$,
and the intersection of this $G$-orbit with $|\cX'|$
equals the $G$-orbit of $|\link_\cX(x)|$.
Since the left-hand vertical map is an inclusion of $G$-spaces,
we have a homotopy pushout diagram of $G$-spaces.
(See e.g.,  \cite{DS95}, \cite{dwyer98sharp}, and \cite{dwyer-henn01}
for basics on pushouts and homotopy pushouts.)
Note also that $|\sta_\cX(x)|$ is $G_x$-contractible,
e.g., by \ref{quillenobs}: since in view of ($\star$),
the mapping taking elements of $\cX_{\leq x}$ to $x$ extends to
a poset endomorphism $F$ of $\sta_\cX(x)$ with image $\cX_{\geq x}$,
which satisfies $F \geq \Id_{\sta_\cX(x)}$;
and then we can send $\cX_{\geq x}$ to its unique smallest element $x$.
Thus if we can show that $|\link_\cX(x)|$ is also $G_x$-contractible,
then the left-hand vertical map is a $G$-homotopy equivalence;
and it follows that right-hand vertical map
is a $G$-homotopy equivalence from $|\cX'|$ to $|\cX|$,
since we have a homotopy pushout of $G$-spaces.
\end{rec}

\begin{prelemma} \label{lemma2}
Let $G$ be a discrete group.
Suppose $\cC$ is a collection of subgroups,
and let $\cC'$ be the subcollection obtained
by removing the $G$-conjugates of some subgroup $P$. 

\begin{enumerate}

\item
Assume for all subgroups $H \leq P$
that the poset $\link_\cC(P)_{\geq H}$ is contractible. \\
Then the inclusion $|E\bO_{\cC'}| \to |E\bO_\cC|$
is a $G$-homotopy equivalence.

\item
Assume for all $H \leq C_G(P)$ that the poset
$\link_\cC(P)_{\leq C_G(H)}$ is contractible. \\
Then the inclusion $|E\A_{\cC'}| \to |E\A_\cC|$
is a $G$-homotopy equivalence.

\item 
Assume for all $H \leq N_G(P)$ that $\link_\cC(P)^H$ is contractible. \\
Then the inclusion $|{\cC'}| \to |\cC|$  is a $G$-homotopy equivalence.

\end{enumerate}

\end{prelemma}

\begin{proof}
To establish (1), we examine the pushout square in \ref{pushoutrec},
with $\cX = E\bO_\cC$ and $x = (G/P,eP)$,
so that $\cX' = E\bO_{\cC'}$, and $x$ has stabilizer $G_x = P$.
We saw in \ref{pushoutrec}
that we get the needed $G$-homotopy equivalence
if we can show that $|\link_{E\bO_\cC}(G/P,eP)|$ is $P$-contractible.
 But for any $H \leq P$, by ($\dagger$),
 $\link_{E\bO_\cC}(G/P,eP)^H$ is equivalent to $\link_\cC(P)_{\geq H}$,
 the contractibility of which is exactly the hypothesis of (1).

The proof of (2) proceeds via a similar pushout square
using the category $E\A_{\cC}$:
here the object $x$ defined by $i : Q \to G$ with $i(Q) = P \in \cC$
has stabilizer $G_x = C_G(P)$
so we need $C_G(P)$-contractibility of $|\link_{E\bO_\cC}(i)|$.
Now ($\dagger$) reduces us to verifying that for all $H \leq C_G(P)$,
 $\link_\cC(P)_{\leq C_G(H)}$ is contractible,
 which again is just the assumption.

Finally (3) is again similar, since the stated assumption
says exactly that $|\link_\cC(P)|$ is $N_G(P)$-contractible,
where $N_G(P)$ is the stabilizer of the object $P$ of $\cC$.
\end{proof}

\begin{lemma} \label{lemma3}
Let $G$ be a discrete group, and let $\cC' \subset \cC$ be collections
such that $\cC \setminus \cC'$ contains
finitely many $G$-conjugacy classes of subgroups.

\begin{enumerate}

\item
Assume for all $P \in \cC \setminus \cC'$
that $\cC_{>P}$ is contractible. \\
Then $|E\bO_{\cC'}| \to |E\bO_\cC|$ is a $G$-homotopy equivalence.

\item
Assume for all $P \in \cC \setminus \cC'$
that $\cC_{<P}$ is contractible.  \\
Then $|E\A_{\cC'}| \to |E\A_\cC|$ is a $G$-homotopy equivalence.

\item
Assume either that 
$\cC_{>P}$ is $N_G(P)$-contractible
for all $P \in \cC \setminus \cC'$, \\
or that
$\cC_{<P}$ is $N_G(P)$-contractible
for all $P \in \cC \setminus \cC'$. \\
Then $|{\cC'}| \to |\cC|$ is a $G$-homotopy equivalence.

\end{enumerate}

\end{lemma}

\begin{proof}
Assume the hypothesis of (1).
We want to argue that we can successively remove
the $G$-conjugacy classes of subgroups in $\cC \setminus \cC'$
in order of {\em increasing} size.
Assume that $P$ is a minimal subgroup in $\cC \setminus \cC'$,
and let $\cC''$ denote the poset obtained
by removing the $G$-conjugates of $P$ from $\cC$.
Then, for $H \leq P$, we have $\cC_{> P} \subseteq \cC_{\geq H}$;
so using ($\star$) we see that 
$\link_{\cC}(P)_{\geq H} = (\cC_{<P})_{\geq H}\ \star\ \cC_{>P}$.
Since $\cC_{>P}$ is assumed contractible,
so is $\link_{\cC}(P)_{\geq H}$.
Hence Lemma~\ref{lemma2}(1) gives
a $G$-homotopy equivalence $|E\bO_{\cC''}| \to |E\bO_{\cC}|$.
Since for all $Q \in \cC'' \setminus \cC'$
we have $\cC''_{>Q} = \cC_{>Q}$ by minimal choice of $P$,
we can continue by induction.

For (2) note that if $P \in \cC \setminus \cC'$,
and $H \leq C_G(P)$, this time $\cC_{< P} \subseteq \cC_{\leq C_G(H)}$,
so using ($\star$) we get
$  \link_\cC(P)_{\leq C_G(H)}
 = \cC_{<P}\ \star\ (\cC_{>P})_{\leq C_G(H)}$;
so that we may remove $G$-conjugates of $P$ by Lemma~\ref{lemma2}(2)
since $\cC_{<P}$ is assumed contractible.
By successively removing conjugacy classes of subgroups
in $\cC \setminus \cC'$
in order of {\em decreasing} size in $\cC \setminus \cC'$,
we conclude that $|E\A_{\cC'}| \to |E\A_\cC|$
is a $G$-homotopy equivalence.

Finally (3) is classical:
First if $\cC_{>P}$ or $\cC_{<P}$ is $N_G(P)$-contractible,
then so is $\link_\cC(P)$.
Thus we may apply Lemma~\ref{lemma2}$(3)$ inductively---removing 
subgroups either bottom-up or top-down as above,
depending on the assumptions on $\cC$.
\end{proof}

\begin{rem}
Note that Lemma~\ref{lemma3} and its proof help explain why 
normalizer formulas tend to have the freedom of choice of $\cC$
of {\em both} the subgroup and the centralizer formulas.
\end{rem}

The next lemma will enable us
to get the vertical lines in Theorem~\ref{equivalencethm}.

\begin{lemma} \label{lemma4}
Let $\cC$ be a collection of $p$-subgroups in a finite group $G$,
with $S \in \rm{Syl}_p(G)$.
\begin{enumerate}

\item
Suppose that $\cC$ is closed under passage to $p$-overgroups. \\
Then the canonical functor $E\bO_\cC \to \cC$
induces an $S$-homotopy equivalence on nerves.

\item
Suppose that $\cC$ is closed under passage to nontrivial subgroups.\\
Then the canonical functor $E\A_\cC \to \cC$
induces an $S$-homotopy equivalence on nerves. 
\end{enumerate}

\end{lemma}

\begin{proof}
Let $H$ denote an arbitrary subgroup of $S$.
By ($\dagger$) the map $E\bO_\cC^H \to \cC^H$ identifies
with the inclusion $\cC_{\geq H} \to \cC^H$.
However if $Q$ and $H$ are $p$-groups and $H \leq N_G(Q)$,
then $QH$ is a $p$-group also normalized by $H$,
so that $QH \in \cC^H$ using the closure hypothesis of (1);
then $Q \mapsto QH$ defines a poset endomorphism
on $\cC^H$ with image in $\cC_{\geq H}$.
Now \ref{quillenobs} shows
that $|\cC_{\geq H}| \to |\cC^H|$ is a homotopy equivalence.
Then (1) follows.

For (2), note that by ($\dagger$) the map $E\A_\cC^H \to \cC^H$
identifies with the inclusion $\cC_{\leq C_G(H)} \to \cC^H$.
By elementary group theory, if $Q$ is non-trivial and $Q \leq N_G(H)$,
then $C_Q(H)$ is non-trivial as well, 
since both $H$ and $Q$ are $p$-groups.
Then $C_Q(H) \in \cC^H$ using the closure hypothesis in (2).
Hence $Q \mapsto C_Q(H)$
gives a poset endomorphism on $\cC^H$ with image in $\cC_{\leq C_G(H)}$.
By \ref{quillenobs}, this induces a deformation retraction
of $|\cC^H|$ onto $|\cC_{\leq C_G(H)}|$.
So $|E\A_\cC| \to |\cC|$ is an $S$-homotopy equivalence as wanted.
\end{proof}

\section{Proofs of the two theorems}   \label{sec:proofthms}

\begin{proof}[Proof of Theorem~\ref{equivalencethm}]
We first establish the horizontal lines.
Since vertical dotted lines always exist (by ($\dagger$) for $H =1$)
it is enough to establish the solid horizontal lines.

Consider the lines between columns $\cB$ and $\cS = \cS_p(G)$,
and between columns $\cBCe$ and $\cCe$.
Note that the $p$-centric condition is closed under $p$-overgroups,
so for $\cX = \cS$ or $\cCe$,
we have $\cX_{>P} = \cS_p(G)_{>P}$ for any $P \in \cX$.
Observe that by elementary group theory, if $Q > P$ then $N_Q(P) > P$;
i.e., $N_Q(P) \in \cS_p(G)_{>P}$.
If $P \not \in \cB$ , then $O_p \bigl( N_G(P) \bigr) \in \cS_p(G)_{>P}$.
Hence the standard inequalities
$$ Q \geq N_Q(P) \leq N_Q(P) O_p \bigl( N_G(P) \bigr)
                \geq O_p \bigl( N_G(P) \bigr) $$
describe a zig-zag of $G$-equivariant functors which by \ref{quillenobs}
show that $\cX_{>P}$ is $N_G(P)$-contractible
to the point $O_p \bigl( N_G(P) \bigr)$.
The two solid lines between the respective columns
now follow, using (1) and (3) of Lemma~\ref{lemma3}.

Next consider the lines between columns $\cS$ and $\cA$:
For $P \in \cS$, we denote by $\Phi(P)$ the Frattini subgroup of $P$,
the smallest normal subgroup of $P$
such that $P / \Phi(P)$ is elementary abelian.
By elementary group theory (see~\cite[Thm.~5.1.1]{gorenstein68}) $\Phi(P) < P$,
and if $Q < P$ then also $\Phi(P) Q < P$.
Furthermore, if $P \not \in \cA$ then $\Phi(P) \neq 1$.
Hence the standard inequalities $Q \leq \Phi(P) Q \geq \Phi(P)$ show
that $\cS_p(G)_{<P}$ is $N_G(P)$-contractible
using \ref{quillenobs}.
The two solid lines between columns $\cS$ and $\cA$
now follow, using (2) and (3) of Lemma~\ref{lemma3}.

We turn to arguments via the functor method of Lemma~\ref{lemma1}.

Consider the lines between columns $\cI$ and $\cS$:
For $P \in \cS$, define
  $$ F(P) = \bigcap_{P \leq S \in Syl_p(G)}\ S ; $$
and observe that $P \leq F(P)$, so that $F \geq \Id_{\cS}$.
Likewise if $P \leq Q \in \cS$,
the Sylow groups above $Q$ are also above $P$, so that $F(P) \leq F(Q)$.
So the two solid lines between columns $\cI$ and $\cS$
now follow using (1) and (3) of Lemma~\ref{lemma1}.
These lines,
together with the earlier solid lines between $\cB$ and $\cS$,
now by composition imply the lines between $\cB$ and $\cI$.

Consider the lines between columns $\cA$ and $\cZ$:
For $P \in \cA$, define $F(P) = \Omega_1 O_p Z \bigl( C_G(P) \bigr)$.
We see that $P \leq F(P)$,
and also that $C_G(P) \leq C_G \bigl( F(P) \bigr)$;
while for $P \leq Q \in \cA$, we have $Q \leq C_G(Q) \leq C_G(P)$,
so that $F(P)$ centralizes $C_G(Q)$ and in particular lies in $C_G(Q)$,
and hence $F(P) \leq F(Q)$.
Thus again $F$ defines an $G$-equivariant poset endomorphism on $\cA$
with $F \geq \Id_\cA$.
We may let $F^{\infty}$ denote the repeated iteration of $F$;
for any finite group $G$, a finite number of iterations suffices.
Then $F^{\infty}$ is idempotent, with image $\cZ$.
The three horizontal solid lines between columns $\cA$ and $\cZ$
now follow using (1)--(3) of Lemma~\ref{lemma1}.
This completes the proof of the horizontal solid lines between columns,
and hence as we mentioned of the dotted horizontal lines as well.

\medskip

We turn to the vertical lines.
We observed
that we at least have dotted vertical lines in all columns.
So it remains to establish the stronger dashed vertical lines:
namely $S$-equivalences for $S$ Sylow in $G$.

We have observed that $\cCe$ is closed under $p$-overgroups,
while the definition of $\cE$ shows it is closed
under nontrivial subgroups;
and by its definition $\cS$ is of course closed both above and below.
Hence the dashed lines in columns $\cCe$, $\cS$, and $\cE$
follow respectively
from (1), (1) and (2), and (2) of Lemma~\ref{lemma4}.
The remaining dashed lines
now follow from the composition
of the parallel vertical dashed lines in the above columns
with adjacent horizontal solid lines already established.
This completes the proof of all the vertical lines shown,
and hence of all the lines in the table of Theorem \ref{equivalencethm}.

\medskip

It remains to show that there are no more lines than those stated.
We provide counterexamples below, considering lines of each type in turn.

{\em Nonexistence of further dotted lines:}
The horizontal and vertical lines established so far
show that we have at most four homotopy types
in the table, represented by $\calD$, $\cCe$, $\cS$, and $\cE$.
To rule out further dotted lines, we will give counterexamples 
showing that these four homotopy types can be distinct.

First consider the group $G = ((\Z/p)^p \times \Z/q) \semi \Z/p$,
for a prime $q$ such that $p \mid q-1$,
where $\Z/p$ acts on $(\Z/p)^p$ via permutation
and on $\Z/q$ by multiplying by a $p$th root unity in $\Z/q$.
This example shows that $\calD$ has a homotopy type
distinct from that of the others,
since in that case $\calD$ is non-contractible,
while $\cCe$, $\cS$, and $\cE$ are contractible.

Next consider $G = \Z/p \times (\Z/q \semi \Z/p)$,
with $q$ as above and the same action.
Here $\cCe$ is non-contractible,
while $\cS$ and $\cE$ (which here equals $\cS$) are contractible.

Finally we claim that $\cS$ is connected while $\cE$ is disconnected
for $G = \SL_2(\F_{p^2}) \semi \Z/p$,
where $\Z/p$ acts by field automorphisms:
Since the inclusion $\SL_2(\F_{p^2}) < G$ gives an identification
of $\cS_p(\SL_2(\F_{p^2}))$ with $\cE_p(G)$,
we see that $\cE$ is disconnected.
To see that $\cS$ is connected,
first observe that the subgroups $S$ and $S'$
of strictly upper and lower triangular matrices in $\SL_2(\F_{p^2})$
are connected in $\cS_p(G)$
(we have a zig-zag $S \leq SQ \geq Q \leq S'Q \geq S'$,
for $Q$ the order $p$ subgroup of field automorphisms,
since $S$ and $S'$ are normalized by $Q$);
this in fact holds for any two Sylow $p$-subgroups in $\SL_2(\F_{p^2})$,
since a Sylow $p$-subgroup is either
equal to $S$ or an $S$-conjugate of $S'$.
 The result now follows since every $p$-subgroup of $G$
 is contained in a Sylow $p$-subgroup,
 and every Sylow $p$-subgroup of $G$
 contains a Sylow $p$-subgroup of $\SL_2(\F_{p^2})$.
 (Interesting counterexamples also arise
from sporadic simple groups for $p=2$;
e.g., $Co_3$ in \cite{benson94},
and other groups such as $M_{12}$ in \cite{BS03}.)
This completes the proof
that there can be no more dotted lines than the ones displayed.

\smallskip

Since we have now established the nonexistence of further dotted lines,
further dashed or solid lines could only arise
within the four homotopy types already indicated. 
For ease of notation, for the rest of the proof
let $D$ denote the dihedral group of order $8$ if $p=2$,
and the extraspecial $p$-group $p^{1+2}_+$ of order $p^3$ and exponent $p$,
if $p$ is odd.

\smallskip

{\em Nonexistence of further dashed lines:}
We work first within the homotopy type of $\cCe$:
Taking $G = D$, so that $S = G$ and $C_G(S) = Z(G)$ is of order $p$,
we see that $\cBCe_{\leq C_G(S)}$ and $\cCe_{\leq C_G(S)}$ are empty,
while $\cCe^S$ is not;
so by ($\dagger$) there are no corresponding dashed lines.
Further taking $H$ to be a maximal abelian subgroup in $G = S$,
we see that $\cBCe_{\leq C_G(H)}$ is empty
while $\cCe_{\leq C_G(H)}$ is not.

Now consider the homotopy type of $\cE$;
again taking $G = D = S$ shows via ($\dagger$)
that $\cE_{\geq S}$ is empty but $\cE^S$ is not.

Next consider the homotopy type of $\cS$.
Again taking $G = D = S$,
and taking $H$ to be a rank-2 elementary abelian $p$-subgroup
shows that $\cB_{\leq C_G(H)}$ and $\cI_{\leq C_G(H)}$ are empty
while $\cA_{\geq H}$ and $\cA^H$ are not.
However $\cA_{\geq S}$ is empty while $\cA^S$ is not.
Finally for $G = \SL_3(\F_p)$ and $S = D$,
$\cB_{\leq C_G(S)}$ is empty but $\cI_{\leq C_G(S)}$ is not.

We are hence left with the homotopy type of $\calD$.
Note that for $G = D = S$, $\calD_{\leq C_G(S)}$ is empty
while $\calD^G$ and $\calD_{\geq G}$ are not.
This reduces us to showing
that there is no dashed line between $\calD$ and $E\bO_\calD$,
which requires a slightly more elaborate argument:

Fix a prime $p$;
and take some $a \geq 1$ for $p > 3$, but $a \geq 3,2$ when $p = 2,3$---the
reason for this choice of $a$ will emerge later.
Now choose a prime $q$ with $p^a \mid q-1$ but $p^{a+1} \nmid q-1$.
Let $G = \Z/q \semi \GL_p(\F_q)$, where $\GL(\F_q)$ acts on $\Z/q$ by letting the index $p$ subgroup $H$ generated by $\SL_p(\F_q)$ and the scalar matrices act trivially and the quotient act faithfully.

In this case, we claim that
$G$ has exactly two conjugacy classes of principal $p$-radical subgroups,
represented by:
a Sylow $p$-subgroup $S$ of order $p^{pa+1}$;
and another subgroup $Q$ of order $p^{a+2}$ given by the central product
of the cyclic subgroup of order $p^a$ of scalar matrices
with an extraspecial group of order $p^3$
(of exponent $p$ for odd $p$, but quaternion of order 8 when $p=2$).
To see this, note that it suffices to work in the quotient group $\GL_p(\F_q)$,
where the statement follows for example using \cite[4A]{AF90} for $p$ odd,
and the analogous result \cite[2B]{an92} for $p = 2$. (The only possible values for the parameters given in those results
are $V_0 = 0$, $s = 1$, $\alpha = 0$;
and either $(e,\gamma) = (p,0)$ corresponding to $S$,
or $(e,\gamma) = (1,1)$ corresponding to $Q$.) By construction $Q$ is a subgroup of $H$, but since $pa > a+2$ (by the choice of $a$) it is not a Sylow $p$-subgroup of $H$. We can hence pick a $p$-subgroup $P \leq N_H(Q)$ strictly containing $Q$. 
Then the poset $\calD^P$ consists of the $G$-conjugates of $Q$ and $S$
normalized by $P$, which is easily seen to be connected:
if $Q'$ is $G$-conjugate to $Q$ and normalized by $P$,
then for $S'$ a Sylow $p$-subgroup containing $Q'P$,
$Q'$ will be connected to $Q$ via $Q' \leq S' \geq Q$;
and likewise for a Sylow $p$-subgroup $S''$ normalized by $P$,
we have $P \leq S''$ so $Q \leq S''$. On the other hand the poset
 $\calD_{\geq P}$ is disconnected, since the fact that $P$ is contained in $H$ and
the non-triviality of the action of $\GL_p(\F_q)$ on $\Z/q$ ensures that there is more than one $G$-conjugate of $S$ containing $P$.

{\em Nonexistence of further solid lines:}
We can continue to work within the indicated homotopy types.
Since there is at least a dotted line in each row of these,
and we have just seen in particular that dotted lines
cannot be strengthened to dashed lines,
it is now enough to see that there can be no new solid lines
going between rows.
For $G = \Z/p \semi \Z/q$, $q \mid p-1$, for $p$ odd or $G = {\rm Alt}_4$ for $p=2$,
and for each $\cC$ in the table,
$\cC_{\leq C_G(G)}$ and $\cC_{\geq G}$ are empty while $\cC^G$ is not.
For $G = \Z/p \times \Z/q$, $q$ a prime different from $p$,
and each $\cC$ in the table,
$\cC_{\geq G}$ is empty, while $\cC_{\leq C_G(G)}$ and $\cC^G$ are not.

This finishes the elimination of any further lines,
and hence the proof of Theorem~\ref{equivalencethm}.
\end{proof}

\begin{rem}\label{histrem1}
The $G$-homotopy equivalences along the normalizer row
in Theorem~\ref{equivalencethm} were at least known classically:
The ordinary equivalence between $\cS$ and $\cA$
was observed by Quillen \cite[2.1]{quillen78},
and between $\cS$ and $\cB$ by Bouc \cite[Cor., p.\ 50]{bouc84}.
The $G$-homotopy equivalence was first observed
by Th\'{e}venaz-Webb \cite{TW91}.
The $G$-equivalence between $\cI$ and $\cS$ was probably first observed
by Alperin sometime in the 1990s via a nerve-of-covering argument;
see also \cite{notbohm01}.
That $\cA$ and $\cZ$ are $G$-homotopy equivalent was observed
in \cite[Sec.~6.6]{benson98}
(though seemingly our argument differs from the one intended there).
Finally, a number of the remaining horizontal equivalences
can be found implicitly in \cite{dwyer98sharp}.
\end{rem}

\begin{proof}[Proof of Theorem~\ref{sharpnessthm}]
The table in the theorem gives references
to the first publication known to us of any particular result.
The two new positive results follow
from the earlier positive results in the same columns
via dashed lines in Theorem~\ref{equivalencethm}.

A counterexample to centralizer sharpness
for $\calD, \cBCe, \cCe, \cB,$ and $\cI$
is provided by $G = D = p^{1+2}_+$,
the extraspecial group of order $p^3$ and exponent $p$ for $p$ odd,
and $D = D_8$ for $p=2$:
This is easy for all the indicated collections $\cC$ except $\cCe$,
since in those other cases $\cC$ consists of just $D$,
and the mod $p$ cohomology of $D$ is different
from that of $Z(D) = \Z/p$.
For $\cC = \cCe_p(D)$ we calculate directly
that $H^*(D;\F_p) \to \lim^0_{P \in \A_\cC} H^* \bigl( C_G(P);\F_p \bigr)$
is not an isomorphism, for instance by observing that the two sides do not 
have the same Krull dimension.

That each of $\cC = \cA$, $\cZ$, and $\cE$ is not subgroup sharp
is likewise easy, since taking $G = \Z/p^2$,
we observe that $H^*(\Z/p^2;\F_p) \to H^*(p\Z/p^2;\F_p)$
is not an isomorphism.
\end{proof}

\begin{rem}
Note that by propagation in the graph of Theorem~\ref{equivalencethm},
the normalizer sharpness of $\cS$ (for example) 
implies 10 other sharpness results.
(See e.g., \cite[\S 2.5]{webb91}, \cite[\S V.3]{AM94},
\cite[7.2]{dwyer98sharp}, or \cite[8.2, \S 9]{grodal02} 
for various proofs of the normalizer sharpness of $\cS$.)
Also note that Theorem~\ref{equivalencethm} allows us for example 
to obtain normalizer sharpness of $\cCe$
from subgroup sharpness of $\cCe$, and vice versa.
(Compare \cite{dwyer98sharp}, \cite{SY97} and \cite{grodal02}.)
\end{rem}

\begin{rem}
We have seen that $|E\bO_\calD|$ and $|\calD|$
are in general not $S$-homotopy equivalent.
However it is possible to analyze the map $E\bO_\calD \to \calD$
to show that it induces an equivalence on Bredon cohomology
with values in any cohomological Mackey functors $\fF$
with $|G:S|$ invertible in $R$
and with the {\em additional} property that
$\fF$ vanish on $p'$-subgroups. 
\end{rem}

\bibliographystyle{amsalpha}
\bibliography{../pcg/poddclassification}

\providecommand{\bysame}{\leavevmode\hbox to3em{\hrulefill}\thinspace}
\providecommand{\MR}{\relax\ifhmode\unskip\space\fi MR }
% \MRhref is called by the amsart/book/proc definition of \MR.
\providecommand{\MRhref}[2]{%
  \href{http://www.ams.org/mathscinet-getitem?mr=#1}{#2}
}
\providecommand{\href}[2]{#2}
\begin{thebibliography}{Web91}

\bibitem[AF90]{AF90}
J.~L. Alperin and P.~Fong, \emph{Weights for symmetric and general linear
  groups}, J. Algebra \textbf{131} (1990), no.~1, 2--22. \MR{MR1054996
  (91h:20014)}

\bibitem[AM94]{AM94}
A.~Adem and R.~J. Milgram, \emph{Cohomology of finite groups}, Grundlehren der
  Mathematischen Wissenschaften [Fundamental Principles of Mathematical
  Sciences], vol. 309, Springer-Verlag, Berlin, 1994. \MR{96f:20082}

\bibitem[An92]{an92}
Jian~Bei An, \emph{{$2$}-weights for general linear groups}, J. Algebra
  \textbf{149} (1992), no.~2, 500--527. \MR{MR1172443 (93j:20025)}

\bibitem[Ben94]{benson94}
D.~Benson, \emph{Conway's group {${\rm Co}\sb 3$} and the {D}ickson
  invariants}, Manuscripta Math. \textbf{85} (1994), no.~2, 177--193.
  \MR{95h:55018}

\bibitem[Ben98]{benson98}
D.~J. Benson, \emph{Representations and cohomology. {II}}, second ed.,
  Cambridge Studies in Advanced Mathematics, vol.~31, Cambridge University
  Press, Cambridge, 1998, Cohomology of groups and modules. \MR{99f:20001b}

\bibitem[Bou84]{bouc84}
S.~Bouc, \emph{Homologie de certains ensembles ordonn\'es}, C. R. Acad. Sci.
  Paris S\'er. I Math. \textbf{299} (1984), no.~2, 49--52. \MR{85k:20150}

\bibitem[Bro75]{brown75}
K.~S. Brown, \emph{Euler characteristics of groups: the {$p$}-fractional part},
  Invent. Math. \textbf{29} (1975), no.~1, 1--5. \MR{52 \#5878}

\bibitem[BS]{BS03}
D.~J. Benson and S.~D. Smith, \emph{Decomposition of {$2$}-completed
  classifying spaces for sporadic simple groups}, in preparation.

\bibitem[DH01]{dwyer-henn01}
W.~G. Dwyer and H.-W. Henn, \emph{Homotopy theoretic methods in group
  cohomology}, Advanced Courses in Mathematics---CRM Barcelona, Birkh\"auser
  Verlag, Basel, 2001. \MR{2003h:20093}

\bibitem[Dre75]{dress75}
A.~W.~M. Dress, \emph{Induction and structure theorems for orthogonal
  representations of finite groups}, Ann. of Math. (2) \textbf{102} (1975),
  no.~2, 291--325. \MR{52 \#8235}

\bibitem[DS95]{DS95}
W.~G. Dwyer and J.~Spali{\'n}ski, \emph{Homotopy theories and model
  categories}, Handbook of algebraic topology, North-Holland, Amsterdam, 1995,
  pp.~73--126. \MR{96h:55014}

\bibitem[Dwy97]{dwyer97}
W.~G. Dwyer, \emph{Homology decompositions for classifying spaces of finite
  groups}, Topology \textbf{36} (1997), no.~4, 783--804. \MR{97m:55016}

\bibitem[Dwy98]{dwyer98sharp}
\bysame, \emph{Sharp homology decompositions for classifying spaces of finite
  groups}, Group representations: cohomology, group actions and topology
  (Seattle, WA, 1996), Proc. Sympos. Pure Math., vol.~63, Amer. Math. Soc.,
  Providence, RI, 1998, pp.~197--220. \MR{99b:55033}

\bibitem[Gor68]{gorenstein68}
Daniel Gorenstein, \emph{Finite groups}, Harper \& Row Publishers, New York,
  1968. \MR{MR0231903 (38 \#229)}

\bibitem[Gro02]{grodal02}
J.~Grodal, \emph{Higher limits via subgroup complexes}, Ann. of Math. (2)
  \textbf{155} (2002), no.~2, 405--457. \MR{1 906 592}

\bibitem[JM92]{JM92}
S.~Jackowski and J.~McClure, \emph{Homotopy decomposition of classifying spaces
  via elementary abelian subgroups}, Topology \textbf{31} (1992), no.~1,
  113--132. \MR{92k:55026}

\bibitem[Not01]{notbohm01}
D.~Notbohm, \emph{Homology decompositions for classifying spaces of finite
  groups associated to modular representations}, J. London Math. Soc. (2)
  \textbf{64} (2001), no.~2, 472--488. \MR{2003a:55022}

\bibitem[Qui78]{quillen78}
D.~Quillen, \emph{Homotopy properties of the poset of nontrivial
  {$p$}-subgroups of a group}, Adv. in Math. \textbf{28} (1978), no.~2,
  101--128. \MR{80k:20049}

\bibitem[SY97]{SY97}
S.~D. Smith and S.~Yoshiara, \emph{Some homotopy equivalences for sporadic
  geometries}, J. Algebra \textbf{192} (1997), no.~1, 326--379. \MR{98e:20024}

\bibitem[TW91]{TW91}
J.~Th{\'e}venaz and P.~J. Webb, \emph{Homotopy equivalence of posets with a
  group action}, J. Combin. Theory Ser. A \textbf{56} (1991), no.~2, 173--181.
  \MR{92k:20049}

\bibitem[Web91]{webb91}
P.~J. Webb, \emph{A split exact sequence of {M}ackey functors}, Comment. Math.
  Helv. \textbf{66} (1991), no.~1, 34--69. \MR{92c:20095}

\bibitem[Yos83]{yoshida83}
T.~Yoshida, \emph{On {$G$}-functors. {II}. {H}ecke operators and
  {$G$}-functors}, J. Math. Soc. Japan \textbf{35} (1983), no.~1, 179--190.
  \MR{84b:20010}

\end{thebibliography}

\end{document}